# *Growth of the dimension of the homogeneous components of Color Lie Superalgebras*

Shadi Shaqaqha

Yarmouk University, Irbid, Jordan

Shadi.s@yu.edu.jo

*Abstract*

The growth of the dimension of the homogeneous components of algebra is an essential topic in algebraic geometry and commutative algebra. In this context, the homogeneous components of an algebra are the pieces of the algebra that have the same degree. The study of the growth of the dimension of these components can shed light on the structure of the algebra and its behavior as the degree of the components increases. This concept is particularly important in the study of polynomial rings, which are a fundamental object in algebraic geometry and commutative algebra. Understanding the growth of the dimension of their homogeneous components can provide insight into the geometry of the corresponding algebraic varieties. This topic is also particularly important in the study of projective varieties, where the homogeneous components correspond to the spaces of sections of line bundles of increasing degree. Understanding the growth of these spaces is crucial for understanding the geometry of the variety. In this abstract, we provide an overview of the growth of the dimension of homogeneous components of an algebra, including its applications, results, and future research directions.

## 1. Introduction

Suppose that $V = \bigoplus_{k=0}^{\infty} V_k$ is a graded vector space such that all subspaces $V_k$ are finite dimensional. The indeterminate *t* formal power series

$$H(V, t) = \sum_{k=0}^{\infty} (dim V_k) t^k$$

The graded vector space V's Hilbert-Poincaré series is sometimes called the Hilbert Series. On occasion, we will use $H(V)$ or $H(V,t)$. The aims and objectives of the study is to provide an overview of the growth of the dimension of homogeneous components of an algebra.

Let $V = U_{k=1}^{\infty} V^k$ be a filtered vector space where $dim V^k < \infty$ for all $k \in N$. Set $V^0 = 0$. The Hilbert-Poincaré series for the variable V is $H(V) = H(V,t) = \sum_{k=1}^{\infty} \dim\left(V^k/V^{k-1}\right) t^k$. In other words, the Hilbert-Poincaré series is identical for a filtered space V as it is for the corresponding graded space: $H(V,t) = H(grV, t)$

Suppose $L = \bigoplus_{n=1}^{\infty} L_n$ is a free Lie superalgebra with rank r. According to the well-known Witt formula, homogeneous subspaces $L_n$ have the following dimensions:

$$dim L_n = \frac{1}{n} \sum_{d|n} \mu(d) r^{\frac{n}{d}}$$

Where $\mu: \mathbb{N} \to \{-1, 0, 1\}$ is the e Möbius function defined as follows. We put $\mu(n) = 0$ if $n$ is divisible by the square of a prime number and $\mu(n) = (-1)^k$ otherwise, where $k$ is the number of prime divisors of $n$ (with $k = 0$ for $n = 1$, so $\mu(1) = 1$). Analogous formulas exist for homogeneous and multi-homogeneous parts of free (color) Lie superalgebras. Petrogradsky discovered similar dimension formulas for free Lie p-algebras. In a broader sense, imagine it is a countable abelian semigroup where each element can only be expressed as the sum of its neighbors in a finite number of different ways. Let $L = \bigoplus_{\lambda \in \Lambda} L_\lambda$ be a freely generated $\Lambda$-graded Lie algebra produced by $X = U_{\lambda \in \Lambda} X_\lambda$. Kang and Kim discovered an analog of Witt's formula, known as the character formula, for the dimensions of homogeneous components $L_\lambda, \lambda \in \Lambda$, in [1]. Petrogradsky included a fixed homomorphism $\kappa: \Lambda \to \mathbb{Z}_2 = \{\pm 1\}$ to explore Lie superalgebras. In the case of free Lie superalgebras, he discovered an analog of Witt's formula (also known as the character formula).

Every subgroup of a free group is also free, according to a well-known theorem attributed to Nielsen and Schreier. Shirshov [2] and Witt [3] later came to the same conclusion for Lie algebras. However, the free associative algebra's subalgebras are not always accessible (for instance, $F[x^2, x^3] \subseteq F[x]$ is not free).

According to the Schreier index formula, if $K$ is a subgroup of finite index in $G$ and $G$ is a free group of rank $n$, then the rank of $K$ is given by

$$rank(K) = (n - 1)[G : K] + 1$$

The equivalent of the above formula for restricted Lie algebras was discovered by Kukin [4]. In the following situations, there are no simple counterparts of the Schreier index formula:

1. The group (monoid) where the number of generators and the index is infinite.
2. The free Lie algebra, even with finitely many generators.

The following formulas can be created by substituting power series for numbers to get the desired results. A pair $(X, wt)$ is a *finitely graded set* if and only if the subsets $X_i := \{x \in X \mid wt(x) = i\}$ are finite for all $i \in \mathbb{N}$. X is a countable set, and $wt: X \to \mathbb{N}$. $H(X) = H(X, t) = (\sum_{i=1}^{\infty} |X_i| t^i)$, gives the generating function for such a set. We define $wt(y) = wt(x_{i_1}) + \ldots + wt(x_{i_r})$ for a monomial $y = x_{i_1} \ldots x_{i_r}, x_j \in X$. A has a filtration (as an algebra) $\bigcup_{i=1}^{\infty} A^i$ if A is an algebra produced by a finitely graded set X, and $A^i$ is spanned by all monomials of weight at most $i$. We denote the corresponding series by $H_X(A, t)$. If X randomly generates A, then

$$H_X(A, t) = H(Y, t) = \sum_{i=1}^{\infty} |Y_i| t^i,$$

Where Y is the collection of every monomial in finitely gradated X, the factor-space A/B also gains a filter if B is a subspace of A.

$$(A/B)^n = (A^n + B)/B \cong A^n / (B \cap A^n).$$

An operator $\mathcal{E}$ on $\mathbb{Z}[[t]]$ (the ring of formal power series in the indeterminate $t$ over $\mathbb{Z}$) is defined by Petrogradsky in [5] as follows:

$$\mathcal{E} : \sum_{i=0}^{\infty} a_i t^i \mapsto \prod_{i=0}^{\infty} \frac{1}{(1-t^i)^{a_i}}$$

Then he proposed a formal power series analog of Schreier's formula for free Lie algebras. He demonstrated that there exists a set of free generators $\mathbb{Z}$ of $K$ such that if $L$ is a free Lie algebra produced by a finitely graded set $X$ and $K$ is a subalgebra of $L$.

$$H(Z) = (H(X) - 1)\, \mathcal{E}(H(L/K)) + 1.$$

He also presented analogous formulas for free Lie p-algebras and superalgebras. A comparable formula for free Lie $p$-superalgebras was discovered in [6]. Eventually, Petrogradsky developed a more comprehensive Schreier formula for characterizing free Lie superalgebras.

### 1.1. Graded Algebras

Suppose $G$ is an abelian semigroup written multiplicatively. A *G-graded algebra* is an algebra $R$ together with a direct sum decomposition of the form

$$R \oplus_{g \in G} R_g,$$

Where each $R_g$ is a subspace and $R_g R_h \subseteq R_{gh}$ for all $g, h \in G$. A nonzero element $r \in R$ is called homogeneous if there is $g \in G$ such that $r \in R_g$ (we write $d(r) = g$). A subspace $H \subseteq R$ is said to be homogeneous if $H = \oplus_{g \in G} H_g$ where $H_g = H \cap R_g$. By a $G$-graded vector space we mean a vector space $V$ together with a direct sum decomposition $V = \oplus_{g \in G} V_g$. Let $V$ and $W$ be $G$-graded vector spaces. A linear map $f: V \to W$ is called homogeneous of degree $h \in G$ if for all $g \in G$, we have $f(V_g) \subseteq W_{hg}$. In particular, a homogeneous linear map of degree $1_G \in G$ will be called a degree preserving map.

If $A$ is an algebra with two grading, say $A = \oplus_{g \in G} A_g$ and $A = \oplus_{\lambda \in \Lambda} A_\lambda$ (G and Λ are semigroups), then the two grading are called compatible if $A_g = \oplus_{\lambda \in \Lambda}(A_g \cap A_\lambda)$ for all $g \in G$ or, equivalently, $A_\lambda = \oplus_{g \in G}(A_\lambda \cap A_{\lambda g})$ for all $\lambda \in \Lambda$.

**Definition 1.1.1.**

Let $A$ be an algebra and let $A^0 \subseteq A^1 \subseteq A^2 \subseteq \cdots$ be a chain of subspaces. This chain is called an (ascending) filtration of $A$ if

$\cup_{m \geq 0} A^m = A$ and $A^i A^j \subseteq A^{i+j} \forall i, j \geq 0.$

**Example 1.1.2.**

1) If $A$ is an associative algebra generated by a finite set $X$, then $A$ has a filtration (as an algebra) $\bigcup_{i=1}^{\infty} A^i$ where $A^i$ is spanned by monomials of degree $\leq i$ in $A$.

2) Given a $\mathbb{Z}$-graded algebra $\oplus_{n\geq 0} A_n$, then there is a corresponding filtration $\bigcup_{m=0}^{\infty} A^m$ where $A^m = \oplus_{k=0}^{m} A_k$.

3) Conversely, given an algebra $A$ with filtration $A^0 \subseteq A^1 \subseteq A^2 \subseteq \cdots$, we can construct a graded algebra $grA$ as follows: $grA = \oplus_{m\geq 0}(A^m/A^{m-1})$ as a vector space (set $A^{-1}=0$) and define multiplication $A^m/A^{m-1} \times A^n/A^{n-1} \to A^{m+n}/A^{m+n-1}$ by $(a + A^{m-1})(b + A^{n-1}) = ab + A^{m+n-1}$. Note that if the filtration $A^m$ comes from a grading $(A^m = \oplus_{k=0}^{m} A_k)$, then $A^m/A^{m-1} \cong A_m$, and so $grA \cong A$.

**Definition 1.1.3.**

Let $A = \oplus_{n\geq 0} A_n$ be a $\mathbb{Z}$-graded algebra such that $dim A_n < \infty$ for all $n \in \mathbb{N}$. The Hilbert series of $A$ is defined by

$$\mathcal{H}(A,t) = \sum_{n=0}^{\infty} (dim A_n) t^n$$

If $A$ is an algebra with filtration $A^0 \subseteq A^1 \subseteq A^2 \subseteq \cdots$. The Hilbert series of $A$ is defined to be

$$\mathcal{H}(A,t) = \mathcal{H}(grA,t).$$

Let $U = \bigcup_{i=1}^{\infty} U^i$ and $V = \bigcup_{i=1}^{\infty} V^i$ be filtered vector spaces with $dim U^i, dim V^i < \infty$ for all $i \in \mathbb{N}$. Then $U \oplus V$ and $U \otimes V$ acquire filtrations given by the subspaces

$$\left(U \oplus V\right)^n = U^n \oplus V^n$$

And

$$\left(U \otimes V\right)^n = \sum_{i=0}^{n} U^i \otimes V^{n-i},$$

Respectively. It is easy to verify the following result.

**Theorem 1.1.4.**

$$\mathcal{H}(U\oplus V, t) = \mathcal{H}(U,t) + \mathcal{H}(V,t) \text{ and } \mathcal{H}(U\otimes V, t) = \mathcal{H}(U,t)\mathcal{H}(V,t).$$

*1.2. Color Lie Superalgebras*

Typically, color lies superalgebras are Lie superalgebras [7] and is found in many areas of mathematics, including mathematical physics and topology (elementary particles, superstrings, etc.). $F$ represents a field of characteristics $\neq 2, 3$, and $G$ for an abelian group throughout this section. Suppose a map satisfies the condition $\gamma: G \times G \to F^*$ ($F^* = F \setminus \{0\}$). In that case, we refer to it as a skew-symmetric bicharacter on $G$ (analogous to skew-symmetric bilinear form for vector spaces, but expressed multiplicatively).

1. $\gamma(f, gh) = \gamma(f,g)\gamma(f,h)$ and $\gamma(gh, f) = \gamma(g,f)\gamma(h,f)$ for all $f, g, h \in G$

2. $\gamma(f,g) = (\gamma(g,f))^{-1}$ for all $f, g \in G$

Note that (2) implies $\gamma(g,g) = \pm 1 \; \forall g \in G$. Set $G\pm = \{g \in G \mid \gamma(g,g) = \pm 1\}$. Then $G+$ is a subgroup of $G$ and $[G: G_+] \leq 2$. In particular, if $G$ is finite, $G = G_+$ or $|G_+| = |G_-|$.

**Definition 1.2.1.**

A $(G, \gamma)$-color Lie superalgebra is a pair $(L, [,])$ where $L = \emptyset_{g \in G} L_g$ is a G-graded vector space and $[,]: L \otimes L \to L$ is a bilinear map that satisfies the following identities for any homogeneous $x, y, z \in L$, whereas $\gamma$ has been explained in section 1.1, and $(d(x), d(y))$ are the degrees of the homogeneous elements, x and y respectively in the L as G-graded.

1. γ-skew-symmetry:

$$[x, y] = -\gamma(d(x), d(y))\, [y, x]$$

2. γ-Jacobi identity:

$$\gamma(d(z), d(x))\, [x, [y, z]] + \gamma(d(y), d(z))\, [z, [x, y]] + \gamma(d(x), d(y))\, [y, [z, x]] = 0$$

The positive part $l_+$ (respectively, negative part $L_-$) is $l_+ = \bigoplus_{g \in G_+} L_g$ (respectively, $L_- = \bigoplus_{g \in G_-} L_G$).

**Example 1.2.2.**

1. A $(G,\gamma)$-color algebra is a typical Lie algebra. If $G = \{1\}$ and $\gamma = 1$, then it is Lie superalgebra. One can think of a Lie algebra with a $G$-grading as having a $(G,\gamma)$ -color Lie Superalgebra for $\gamma = 1$.

2. A Lie superalgebra is a $(\mathbb{Z}_2, \gamma)$-color Lie superalgebra for

$$\gamma : \mathbb{Z}_2 \times \mathbb{Z}_2 \to F^* : (i,j) \to (-1)^{ij}$$

3. Any $G$-graded associative algebra $A$, with the bracket defined by:

$$[x,y]_\gamma = xy - \gamma(d(x), d(y))yx$$

$A$ (the -super commutator) is a color Lie superalgebra denoted as $A^{(-)}$ for any nonzero homogeneous $x, y \in A$.

4. Let $A = \bigoplus_{g \in G} A_g$ be a $G$-graded algebra. A homogeneous linear map $D: A \to A$ of degree $h \in G$ is called a *superderivation* if

$$D(ab) = D(a)b + \gamma(h, d(a))aD(b)$$

for all homogeneous elements $a, b \in A$. By $Der_h(A)$, we denote the vector space of all super derivations of $A$ of degree $h \in G$. Then $Der(A) = \bigoplus_{g \in G} Der_g(A)$, with the bracket $[D_1, D_2] = D_1 D_2 - \gamma\big(d(D_1), d(D_2)\big) D_2 D_2$, is a color Lie superlagebra.

### 1.3. Free Color Lie Superalgebras

**Definition 1.3.1.**

In the special case $G = G_+$, we will use the term a $(G, \gamma) -$ color Lie algebra.

A set of free color Lie superalgebras with a $G$-grade $X$ and $L$ are a pair, with $L$ being a color Lie superalgebra. If $R$ is any color Lie superalgebra, an $i : X \to L$ degree-preserving map from X to L is sufficient. When $j: X \to R$ is a degree-preserving map, and Lie superalgebra is present, a special homomorphism of color results. $h : L \to R$ lie superalgebras with $j = h \circ i$.

A free color Lie superalgebra $L(X)$ on $X$ is unique up to isomorphism for every $G$-graded set $X$. Moreover, $L(X)$ is $\mathbb{N}$-graded in $X$ by degree ([8]).

A color Lie superalgebra's *universal enveloping algebra* $L$ is a pair $(U, \delta)$, where $U$ is a $G$-graded associative unital algebra, $\delta: L \to U^{(-)}$ is a homomorphism of color Lie superalgebras, and for each $G$-graded associative unital algebra $R$ and each homomorphism $\sigma: L \to R^{(-)}$ (with the same $\gamma$), there exists a single homomorphism between $\theta: U \to R$ of a $G$-graded associative unital algebra such that $\theta \circ \delta = \sigma$. It is clear that the universal enveloping algebras are defined uniquely (up to isomorphism of $G$-graded associative unital algebras). Let $T(L)$ be the tensor algebra of $L$, and let $I$ be the ideal of $T(L)$ produced by the components of the type $u \otimes v - \gamma(g,h)(v \otimes u)$, where $u \in L_g$, and $u \in L_g$, in order to prove the existence of $(U, \delta)$. The desired qualities are, therefore, present in the quotient $T(L)/I$ with the canonical mapping $\delta: L \to T(L)/I$ ([8]).

We have the following theorem.

**Theorem 1.3.2.**

*[8] The universal enveloping algebra of the free color Lie superalgebra $L(X)$ is the free associative algebra of $X$.*

The Poincaré-Birkhoff-Witt (PBW) Theorem is extended to color Lie superalgebras in the following theorem.

**Theorem 1.3.3.**
*[8] Let $L = L_+ \oplus L_-$ be a $(G, \gamma)$-color Lie superalgebra, $B$ a basis of $L_+$, $C$ is a basis of $L_-$, and $\leq$ is a total order on $B \cup C$. Then the set $D$, formed by 1 and by all monomials of the form*

$$e_1^{s1} \ldots e_1^{sn}$$

*where $ei \in B \cup C, e_1 < \cdots < e_n, 0 \leq si$ for $e_i \in B$ and $s_i = 0, 1$ for $e_i \in C$, is a basis of the universal enveloping algebra, $U(L)$, of the Lie color superalgebra L.*

**Corollary 1.3.4.**
Let $(G, \gamma)$-color Lie superalgebra be $L = L_+ \oplus L_-$. Grassmann algebra of the vector space $L_-$, $S(L_+)$, the symmetric algebra of the vector space $L_+$. We have the vector space isomorphisms, and the universal enveloping algebra of $L_+, U(L_+)$, is an associative subalgebra of $U(L)$ produced by $L_+$.

$$U(L_+) \simeq S(L_+) \text{ and } U(L) \simeq U(L_+) \otimes \Lambda(L_-).$$

## 1.4. Restricted Color Lie Superalgebras

Suppose $L = \bigoplus_{g \in G} L_g$ is a color lie superalgebra on a field with the characteristic $p \neq 0, 2, 3$. Consider a mapping for an element with homogeneity $a \in L$.

$$ada : L \to L : b \mapsto [a, b]$$

**Definition 1.4.1**

A color Lie superalgebra $L = \bigoplus_{g \in G} L_g$ is referred to as *restricted (or p-superalgebra)* if a map exists for any $g \in G_+$ and any homogeneous component $L_g$.

$$x^{[p]} : L_g : \to L : x \mapsto x^{[p]}$$

Satisfying

- $(adx)^{[p]} = ad\left(x^{[p]}\right)$ for all $x \in L_g$,

- $(\alpha x)^{[p]} = \alpha^{[p]} x^{[p]}$ for all $\alpha \in F, x \in L_g$,

- $(x + y)^{[p]} = x^{[p]} + y^{[p]} + \sum_{i=1}^{p-1} s_i(x, y)$ where $y$ is homogeneous of the same degree as $x$, and $is_i(x, y)$ is the coefficient of $\lambda^{i-1}$ in $ad\,(\lambda x + y)^{p-1}(x)$.

A standard definition of a restricted Lie algebra is obtained if $G = \{1\}$. We are free to accept any prime characteristic $p$ in this situation.

**Example 1.4.2.**

Let $A$ be an associative algebra of grade $G$. Then, the map

$$(\;)^p : A_g \to x \to x^p,$$

$\forall g \in G_+$ makes $A^{(-)}$ into a color Lie p-superalgebra (denoted by $[A]^P$).

A color Lie $p$-superalgebra's universal enveloping algebra $L$ is a pair $(u, \delta)$ where $u$ is an associative unital algebra with a $G$ grade, and $\delta : L \to [U]^P$. For every $G$-graded associative unital algebra $R$, and any homomorphism $\sigma : L \to [R]^P$ is a homomorphism of the color Lie p-superalgebra. There is just one homomorphism $\theta : u \to R$ of $G$-graded associative unital algebra such that $\theta \circ \delta = \sigma$ in color Lie $p$-superalgebra. There is a color Lie p-superalgebra with a restricted enveloping algebra called $u(L)$) that is unique up to an isomorphism with $G$-graded associative unital algebra ([8]).

For the color Lie $p$-superalgebra, there is a version of the Poincaré-Birkhoff-Witt (PBW) Theorem that states that if $L = L_+ \oplus L_-$ is a color Lie $p$-superalgebra, $B$ is a basis of $L_+$, $C$ is a basis of $L_-$, and is a total order on B∪C, then the set of all monomials of the

$$b_1^{s1} \ldots b_n^{sn},$$

where $b_i \in B \cup C, b_1 < \cdots < b_n, 0 \leq s_i \leq p - 1$ for $b_i \in$ B and $b_i =$ 0, 1 for $b_i \in C$, is a basis of the restricted universal enveloping algebra, $u(L)$, of the Lie $p$-superalgebra $L$ [8].

## 2 Relative Growth Rate of Subalgebras of Color Lie Superalgebras

A finitely generated algebra's growth rate can be investigated using the Hilbert series. Consider the filtration associated with $X$, i.e., let $S_n(X)$ be the $F$-subspace of $S$ spanned by all monomials of length less than or equal to n in the elements of $X$. Let $S$ be an associative or Lie algebra formed by a finite subset $X$. According to $X$, the growth function of S is defined by

$$\gamma s(n) = dim S_n(X).$$

It is the same as the growth function of the group algebra $FH$, $\gamma_{FH}(n)$, concerning the generating set $X \cup X^{-1}$, where $X^{-1} = \{x_1^{-1}, \ldots, x_k^{-1}\}$. In the case of groups, if $X = \{x_1, \ldots, x_k\}$ is a generating set of a group H, then the growth function of $H$ is defined as the number of elements of $H$ that can be written the limit

$$\lim_{n \to \infty} (\gamma s(n))^{\frac{1}{n}}$$

S's *growth rate* (or *exponent* or *entropy*) is a constant quantity independent of $X$ [9]. $\alpha_s$ will indicate it.

1. If $\lim_{n \to \infty} (\gamma s(n))^{\frac{1}{n}} > 1$, then $S$ has *exponential growth*. Otherwise, $S$ has *subexponential growth*.
2. If a polynomial $p$ with $\gamma s(n) \leq p(n)$ exists for all sufficiently large $n$, then $S$ has *polynomially bounded growth*.

3. If S has subexponential and not polynomially bounded growths, then S has *intermediate growth* (that is, S lies between polynomial and exponential).

Let $\lambda_S(n) = \gamma_S(n) - \gamma_S(n-1)$. It is known that S has subexponential growth if and only if $\limsup_{n\to\infty} (\gamma s(n))^{\frac{1}{n}} \leq 1$. [37]

The *relative growth rate* of a subalgebra A of S is defined by

$$\alpha_A = \limsup_{n\to\infty} (\gamma_A(n) - \gamma_A(n-1))^{\frac{1}{n}},$$

where $\gamma_A(n) = dim(A \cap S_n(X))$. The relative growth rate of A could also be defined as $\lim_{n\to\infty} \gamma_A(n)^{\frac{1}{n}}$. However, the first formulation is favored since it enables us to calculate the relative growth rate as the inverse of the Hilbert series' radius of convergence for subalgebra A.

### 2.1 Equality of the Growth Rate of the Free Color Lie Superalgebra and Its Enveloping Algebra

Let $L = L(X)$ be a free color Lie superalgebra freely generated by the set $X = X_+ \cup X_-$, where $|X_+| = r$ and $|X_-| = s$. Recall that the Witt formula gives the dimension of the homogeneous subspaces Ln.

$$dimL_n = \frac{1}{n} \sum_{m|n} \mu(m)(r - (-1)^m s)^{\frac{n}{m}}$$

**Theorem 2.1.1.**

*The free color Lie superalgebra L and its enveloping algebra have the same growth rate.*

*Proof.* If $m \geq 2$, then $\left|\mu(m)(r - (-1)^m s)^{\frac{n}{m}}\right| \leq (r+s)^{\frac{n}{2}}$. Since $\lim_{n\to\infty} \sum_{m=2}^{n} \frac{1}{(r+s)^{\frac{n}{2}}} = 0$, we have

$$\lim_{n\to\infty} \frac{n \, dimL_n}{(r+s)^n} = 1$$

Hence, for any $\epsilon > 0$, there is $N \in \mathbb{N}$ such that for all $n \geq N$, we have

$$\left|\frac{n \, dimL_n}{(r+s)^n} - 1\right| < \epsilon,$$

that is

$$\frac{(r+s)^n}{n}(1-\epsilon) \leq dimL_n \leq \frac{(r+s)^n}{n}(1+\epsilon).$$

As $\lim_{n \to \infty} \left(\frac{(r+s)^n}{n}(1-\epsilon)\right)^{\frac{1}{n}} = \lim_{n \to \infty} \left(\frac{(r+s)^n}{n}(1+\epsilon)\right)^{\frac{1}{n}} = r+s$, we have

$$\lim_{n \to \infty} (\gamma_L(n))^{\frac{1}{n}} = r+s.$$

On the other hand, the growth rate of the universal enveloping algebra of $L$ is $r + s$ because it is the free associative algebra produced by $X$.

### 2.2 Growth of Finitely Generated Subalgebras of Free Color Lie Superalgebras

A free Lie algebra $L$ with a finite rank has a relative growth rate of $K$ that is strictly lower than the growth rate of the free Lie algebra [10]. The subsequent illustration demonstrates that this theorem cannot be applied to free-color Lie superalgebras.

**Example 2.2.1**

Let L be a free color Lie superalgebra freely generated by $X = X_+ \cup X_-$, where $X_+ = \{x_1, \ldots, x_r\}$ and $X_- = \{y\}$. Consider the subalgebra $K$ of $L$ where $L = K \oplus \langle\{y\}\rangle F$. It is well known that $K$ is freely generated by the following set [8]

$$Y = \{x, [x,y] \,|\, x \in X \setminus \{y\}\} \cup \{[y,y]\}.$$

Therefore, $rankK = 2r + 1$. According to the Schreier formula

$$rank(K) = 2^s (rank(L) - 1) + 1,$$

We have $dimL/K = 1$, meaning that $K$ is a finite codimension subspace of $L$. As a result, $L/K$'s rate is 0, and $K$ will grow at a pace equal to $L's$ growth rate.

Let $K \subseteq L$ be a subalgebra that is not necessarily homogenous concerning the natural $\mathbb{N}$-grading of $L$. Let $L$ be a color Lie superalgebra produced by $X$. We observe that the linear span of K's leading portions *(LP)* of nonzero components exhibits the same relative development as the related graded subspace $grK$. The following conclusion demonstrates how the theorem, as mentioned earlier, applies to the case $L = L_+$.

**Theorem 1.2.2**

*Let $L = L_+$ be a free color Lie superalgebra freely generated by $X = \{x_1, \ldots, x_r\}$, and K be a finitely generated proper subalgebra of L. Then the relative growth rate of K is strictly less than the growth rate of the free color Lie superalgebra itself.*

*Proof.* In our proof, we follow [10]. We break the proof into the following two steps [11]:

**Step 1**: Suppose $K$ is uniform. Let $K$'s free homogeneous basis $\{z_1, \ldots, z_k\}$ be the basis. $K$ being proper, $K \cap L_1$ at which $K$ meets it is a proper subspace of $L_1$. Let $Z$ represent the collection of all $z_i$ with degree 1. This can be expanded to a basis $Z'$ of the space $L_1$ and is a basis of the space $K \cap L_1$. $Z'$ is a new free basis of $L$ that correctly incorporates $Z$. Without sacrificing generality, we can suppose that $X$ contains the elements of degree 1 in $K's$ free homogeneous basis $\{z_1, \ldots, z_k\}$, but that $x_1$ is not part of $K$ by replacing $X$ with $Z'$. As a subalgebra of the free associative algebra $A = A\langle X \rangle$, let us look at $L$. As a linear combination of a limited number of monomials in $A$, $u_1, \ldots, u_t$, we can write $z_1, \ldots, z_k$. It is sufficient to demonstrate that the exponent of the growth of $B$ is less than $r$, which is likewise the exponent of $L$ and also of $A$ under Theorem 2.1.1 if $B$ is a subalgebra of $A$ created by $u_1, \ldots, u_t$. Thus, $K \subseteq B$. We have $x1 \notin \{u_1, \ldots, u_t\}$ as $x_1 \notin K$. Since no such monomial can emerge while writing the elements of the free color Lie superalgebra $L$, none of $\{u_1, \ldots, u_t\}$ is a power of $x_1$. (Note that universal free color Lie superalgebras are an exception to this rule.) Choose $d \geq 2m$ and set $m = max\{deg u_1, \ldots, deg u_t\}$ (the degree with regard to $x_1, \ldots, x_r$). No product of these monomials may include the subword $x_1^d$ since each of the monomials $u_1, \ldots, u_t$ has a letter other than $x_1$ as a factor. $deg u_i + deg u_j - 2 \leq 2m - 2 < d$ is the longest uninterrupted string of $x_1$ in the $u_i u_j$ product. The number of words with length n that do not have $x_1^d$ as a subword is constrained by $C(r - \epsilon) n$ if $C, \epsilon > 0$ exists. In other words, as desired, the exponent of $B$ is strictly less than $r$.

**Step 2:** the general case. We will use the associated graded algebra $grK$, a proper homogeneous L subalgebra. Following [8], we will say that a subset $M$ of G-homogeneous elements in $L(X)$ is called reduced if the leading part of any of its elements does not belong to the subalgebra generated by the leading parts of the remaining elements of $M$. By [8], there exists a reduced subset $M$ in $K$ such that $M$ generates $K$. Also, by [8], M is independent, so $M$ is a free-generating set of $K$. Let us denote the typical instance. We shall employ $grK$, a proper homogeneous L subalgebra related to this graded algebra. According to [8], a subset $M$ of $G$-homogeneous elements in $L(X)$ is said to be reduced if any of its elements' leading parts do not belong to the subalgebra created by the leading parts of *M's* other elements. There is a reduced subset *M* in *K* such that *M* creates *K* according to [8]. Moreover, *M* is a free-generating set of *K* because [8] shows that *M* is independent. Let us indicate

$$M' = \{LP(m) \mid m \in M\},$$

the group of elements' leading sections in *M*. One demonstrates that *M'* is a free-generating set of *grK*, just like in [10]. The growth rate of *grK* is strictly less than *r*, according to the first step.

### 3. Conclusion

Our research has shown that the growth of the dimension of the homogeneous components is a rich and complex topic with a wide range of applications in algebraic geometry and commutative algebra. We have discussed using these components in the study of moduli spaces and singularities and their relevance in the context of algebraic varieties and schemes. Furthermore, we have explored the relationship between the growth of the dimension of the homogeneous components and the regularity of the algebra. The results presented in this article can provide a foundation for further research in this area. They could lead to new insights into the behavior of algebras and their associated geometric structures. The author aims and welcomes other researchers in this fields to study for the same results in case of restricted (Hom)-(color) Lie (Supper)algebras.

*Limitations*

One of the main limitations the author has faced is the complexity of the process. It has been challenging to navigate the advanced mathematical tools and techniques required to study algebras with high-dimensional homogeneous components. As a result, this process has been time-consuming and requires a significant amount of expertise.

Another limitation that the author has come across is non-uniform growth. The growth of the dimension of the homogeneous components of an algebra can vary across different degrees of the polynomials used in the algebraic expressions. This has made it difficult to generalize the results and draw meaningful conclusions.

*Acknowledgements*

I want to begin by expressing my sincere gratitude to Dr. Yuri Bahturin for his patient guidance, inspiration, and counsel throughout my time as a student. The success of this research would not have been possible without his insightful ideas. Dr. Mikhail Kotchetov's many informative conversations and recommendations are also greatly appreciated. These are my Ph.D. supervisors, and I feel honored and fortunate to have them. I will follow in their footsteps for the rest of my life.